%% file: Treisman-Rationally_generated_arcs.tex
\documentclass[11pt]{amsart}
\usepackage{fullpage}
\usepackage[all]{xy}

\input{macros}

\begin{document}

\title{Rationally generated arcs}
\author{Zachary Treisman}
\email{treisman@math.tifr.res.in}
\address{Department of Mathematics\\
  University of Washington\\
  Box 354350\\
  Seattle, Washington 98195\\
  USA}
\curraddr{School of Mathematics\\
  Tata Institute of Fundamental Research\\
  Homi Bhabha Road, Mumbai 400 005\\
  India}
\thanks{This research was done in part while the author was supported by the University of Washington's National Science Foundation VIGRE grant.}
\maketitle
\tableofcontents

\begin{abstract}One distinguishing feature of rational curves is that they have algebraic parameterizations.  Arc spaces are a way of describing approximations to parameterizations of all curves in some fixed space. Playing on these descriptions, this paper introduces {\em rationally generated arcs}, and applies them towards the description and enumeration of rational curves in algebraic varieties.
\end{abstract}

\section{Introduction}

Arc spaces represent approximations of parameterized curves in an algebraic variety, and as such they should tell us a lot about rational curves, when they exist, since these are the only curves with algebraic parameterizations.  For some varieties, there is a nice correspondence between spaces of arcs and rational curves.  Some properties of varieties defined via tangents have refinements in terms of arcs.  And in the  case of $X\subset\Pp^n$ the relationship between arcs and the curves they approximate is simple enough to be utilized for explicit calculations. 

In section \ref{rat_gen} {\em rationally generated arcs} are defined. Roughly speaking, an arc $\gamma$ is rationally generated if there is a rational curve $C$ such that $\gamma$ is an approximation of $C$.  One interesting aspect of this idea is the following chain of implications.

\begin{thm}\label{uni_rga_rc}
If $X$ is any variety, then
$$
{\rm \begin{tabular}{c} $X$ is \\unirational\end{tabular}}\implies{\rm \begin{tabular}{c} $X$ has generically\\ rationally generated arcs\end{tabular}}\implies{\rm \begin{tabular}{c} $X$ is rationally\\ connected\end{tabular}}
$$ 
\end{thm}

Recall that there are no examples known distinguishing rational connectedness from unirationality, but it is conjectured that such spaces exist.  One reason why arcs may be useful in this context is that the explicit descriptions of arc spaces allow for rational generation properties to be checked in many cases.  Specifically, Theorem \ref{locus} calculates the arcs that genereate lines in a complete intersection.  

\begin{notn}
$\widehat\Pp^n_1$, to be defined in \ref{GGarcs}, is the bundle of projectivized arcs of degree 1 on $\Pp^n$. There are natural classes $h_0$ and $h_1$ in the intersection ring $A(\widehat\Pp^n_1)$ corresponding to the hyperplanes in $\Pp^n$ and the relative hyperplanes in the bundle of 1-arcs.  These classes will be defined precisely in section \ref{intersect}. 
\end{notn}

\begin{thm}\label{locus}
Consider $X\in \Pp^n$, a generic complete intersection of type $(d_1,\ldots,d_r)$. The locus of 1-arcs in $\widehat\Pp^n_1$ which are rationally generated by lines in $X$ has cycle class
$$ 
\prod_{j=0}^r\prod_{i=0}^{d_j}\left( ih_0+(d_j-i)h_1\right) \in
A(\widehat\Pp^n_1).
$$ 
In particular, when $\sum d_i=2(n-1)-r$ the class that is the result of
the calculation is
$$
Nh_0^{n-1}h_1^{n-1}.
$$ 
$N$ is computable in terms of $r$, $n$ and the $d_i$.
This class indicates $N$ lines in $\Pp^n$.\end{thm}

Among the examples recovered by this result are the 2875 lines on a quintic 3-fold, and the 6 lines through a general point of a cubic threefold.  This second example is an instance of the following result, which refines the upper bound of $m!$ on the number of lines through a general point of an $m$--dimensional variety published in \cite{Land}.

\begin{thm}\label{genpt} 
Consider $X\subset\Pp^n$, a variety of dimension $m$, such that $X$ has a finite number of lines through a generic point, and is a component of some complete intersection that is general in the sense of section \ref{cigenpt}.
\begin{itemize}
\item 
The number of lines through a general point $x\in X$, if finite, is at
most $m!$.

\item 
If $X$ is not a hypersurface, then this bound improves to $2(m-1)!$.

\item 
If $X$ is not a hypersurface and is not contained in a quadric, then
the bound improves to $3!(m-2)!$

\item 
In general, if $X$ is a variety of codimension $r$, the number of
lines contained in $X$ and passing through a general point $x\in X$ is bounded above by
$$ 
\max\left\{\,\,\prod_{i=1}^rd_i!\,\, \vline
\,\begin{array}{l}\sum_{i=1}^rd_i=N-1,\,d_i>1,\\
                  \mbox{\rm and $X$ is contained in a}\\
                  \mbox{\rm complete intersection}\\
                  \mbox{\rm of type $(d_1,\ldots,d_r)$}\end{array} \right\}.
$$ 
\end{itemize}
\end{thm}

It would be interesting to extend this technique to deal with the description of curves of higher degree.  The difficulty in this aspect of the project consists of finding explicit descriptions of the intersection rings of the relevant arc spaces.

I am working over the field $\C$ of complex numbers throughout.  Much of this machinery does not depend on this restriction; this will hopefully be the subject of some later investigation.

{\bf Acknolwedgements.}  Sincere thanks are due to my PhD advisor, S\'andor Kov\'acs. In addition, I had a number of illuminating and encouraging discussions with Ravi Vakil, Wim Veys, and Carolina Araujo.

\section{Arc Spaces}

The appearance of arc spaces on the modern mathematical scene originates in the early work of John Nash, which originally appeared in 1968 and is published in \cite{Nash}.  See \cite{Blik} and \cite{Veys} for very readable introductions.

Arc spaces are sometimes called jet spaces.  This terminology will be avoided here, as the term jet space or jet bundle has also been used in other contexts. 

\begin{defn}\label{affine_arc_space_1}
Arc spaces are moduli spaces of infinitesimal curves.  Consider a scheme $X$ and the functors $Z\mapsto\Hom(Z\times {\spec\C[\eps]/(\eps^{m+1})},X)$. These functors are representable.  That is, for an integer $m\geq0$,
$$
\Hom(Z\times {\spec\C[\eps]/(\eps^{m+1})},X)=\Hom(Z,X_m).
$$ 
$X_m$ is called the $m^{th}$ arc space of $X$.
\end{defn}

On the level of closed points, 
$$
X_m=\Hom{}(pt,X_m)=\Hom{}(\spec\C[\eps]/(\eps^{m+1}),X),
$$ 
and for any arc $\gamma\in X_m$, $\gamma(0)$ is a closed point of $X$.
 
Note that $X_0\iso X$ and $X_1\iso TX$.  Over smooth points $x\in X$ $X_m$ is a $TX_x$-bundle over $X_{m-1}$.  

\begin{notn}
A lower index on an arc space coordinate indicates the power of $\eps$ of which that coordinate is a coefficient.  The arc space index will always be written last.  

For example, $x_0,\ldots,x_n$ are coordinates on $\Pp^n$, and $x_{00},x_{01},\ldots,x_{n0},x_{n1}$ are coordinates on $\Pp^n_1$, and a point in $\Pp^n_1$ expressed in these coordinates is $[x_{00}+x_{01}\eps:\ldots:x_{n0}+x_{n1}\eps]$.
\end{notn}

\begin{notn} The fiber over a point $x\in X$ will be denoted $X_{m,x}$.
\end{notn}

\begin{ex}\label{jetideal}\label{affine_arc_space_2} 
Suppose $X=\spec\frac{\C[x_1,\ldots,x_n]}{(f_1,\ldots,f_s)}$ is an affine variety in $\A^n$.  The arc space $X_m$ is defined by the equations
$$
f_r(x_1(\eps),\ldots,x_n(\eps))\equiv0 \,\,(\mbox{mod}(\eps^{m+1})) \quad r=1,\ldots,s.
$$ 
In the ring $\C[x_{ij}|i=1,\ldots,n\, j=1,\ldots, m]$, each $f_r$ defines $m+1$ equations, the coefficients of the $\eps^j$ as $j=0,\ldots, m$.  Hence, writing $f_r(x_1(\eps),\ldots,x_n(\eps))=\sum f_{rj}\eps^j$,
$$
X_m=\spec\frac{\C[x_{ij}|i=1,\ldots,n\, j=1,\ldots, m]}{(f_{10},\ldots,f_{sm})}.
$$

If $X$ is an arbitrary quasi-projective variety with affine cover $\{U_\alpha\}$, the local descriptions for $U_{\alpha m}$ patch together to provide equations for $X_m$.
\end{ex}
\label{projarcs}\label{ggarcs}

In the context of this paper, describing curves that exist globally  on $X$, two problems with these arc spaces present themselves.  First, they have too much information, since for many considerations the local parameterizations that a curve may have are less relevant than the locus that is swept out globally when the local data is patched
together.  Second, they are not projective, so lots of powerful global machinery from algebraic geometry is unavailable.

These issues are remedied by descending to projectivized
arc spaces.  The projectivization is accomplished via taking a quotient by reparameterizations of the arcs in the arc space.

\begin{defn}\cite{G-G},\cite{Voj}\label{GGarcs}  With notation as in \ref{jetideal}, define the {\em (Green--Griffiths) projective arc spaces},  
$$
\widehat X_m=\proj\frac{\C[x_{ij}|i=1,\ldots,n\, j=1,\ldots, m]}{(f_{10},\ldots,f_{sm})}.
$$ 
\end{defn}

As above, this local definition serves to define the projective arc spaces for an arbitrary quasi-projective variety or scheme by patching together over an affine cover.  $\widehat X_m$ is a quotient of $X_m$ by the $\C^*$ action $\eps\mapsto\lambda\eps$.

There are other candidate projective arc spaces, such as the Demailly-Semple spaces discussed in \cite{Dem}.  However, different curves can correspond to the same Demailly-Semple arc, so their utility in this present context is less obvious.

\section{Rationally generated arcs}\label{rat_gen}

\begin{defn}A {\em parameterization} of a curve $C\in X$ is a morphism of schemes
$$
\spec \C[[t]]\to X
$$
such that the image is contained in $C$.
\end{defn}

\begin{defn}The {\em Taylor series maps} are induced by truncations of power series:  $\C[[t]]\to\C[\eps]/(\eps^{m+1}); \quad t\mapsto\eps$.
$$
\tau:\left\{\mbox{\rm \begin{tabular}{c}parameterized\\ curves in $X$\end{tabular}}\right\}\longrightarrow\widehat X_m
$$
\end{defn}

\begin{defn} Say that a scheme $X$ has {\em rationally generated $m$-arcs} if the Taylor series map $\tau_x:\Hom(\Pp^1,X;0\mapsto x)\to\widehat X_{m,x}$ is surjective at all points $x\in X$.  Say that $X$ has {\em rationally generated arcs} if $\tau_x$ is surjective for all $x\in X$ for all $m$. \end{defn}

If an ample line bundle $\sL$ is chosen on $X$, one can further refine this condition by requiring that the degree of the rational curve generating an arc be of some bounded degree, measured with respect to this $\sL$.

\begin{defn} A scheme $X$ with an ample line bundle $\sL$ has {\em $d$-rationally generated $m$-arcs} if  $\tau_x:\Hom_d(\Pp^1,X;0\mapsto x)\to\widehat X_{m,x}$ is surjective for all $x\in X$, where $\Hom_d(\Pp^1,X;0\mapsto x)=\{\fee\in\Hom(\Pp^1,X;0\mapsto x)\,|\,\deg\fee\leq d \}$.
\end{defn}

The adverb {\em generically} preceding any of the conditions just defined indicates that the condition is true on some open subset $U\subset X$.

Translated into this language, Taylor's theorem for rational curves in projective space says that $\Pp^n$ has $d$-rationally generated $d$-arcs for all $d$, and that there is in fact a unique degree $d$ parameterization of a rational curve for every arc in $\Pp^n_d$.  

\begin{notn}\label{tilde}
For an arc $\gamma\in\Pp^n_m$ written in coordinates as
$$
\gamma=[x_{00}+\cdots+x_{0m}\eps^m:\cdots:x_{n0}+\cdots+x_{nm}\eps^m]
$$
Define the parameterized rational curve $\theta(\gamma(t))$ via the vector space map $\eps\mapsto t$.
$$
\theta(\gamma(t))=[x_{00}+\cdots+x_{0m}t^m:\cdots:x_{n0}+\cdots+x_{nm}t^m]
$$
And write $\tilde\gamma$ for the curve that is parameterized by $\theta(\gamma(t))$.
\end{notn}

Note that the tilde construction is well defined on the Green-Griffiths projective arc spaces.  That is, if $\gamma, \gamma'\in X_m$ are such that $\bar\gamma=\bar\gamma'$ under the projection $X_m\to\widehat X_m$ then $\tilde\gamma=\tilde\gamma'$.  The same notation $\tilde\gamma$ will be used regardless of whether $\gamma\in X_m$ or $\gamma\in \widehat X_m$.

\begin{ex}\label{tildeline1}
If $X=\Pp^n$ and $\gamma\in\widehat X_1$ with coordinate representation $(\eps,0,\ldots,0)$ on the chart $U_0=\{x_0\neq0\}=\{(\frac{x_1}{x_0},\ldots,\frac{x_n}{x_0})\} $ then $\tilde\gamma$ is the line $(t,0,\ldots,0)$.
\end{ex}

\begin{ex}\label{tildeline2}
Again consider $X=\Pp^n$ and $\gamma\in\widehat X_1$, this time with coordinate representation $(\eps,0,\ldots,0)$ on the chart $U_1=\{x_1\neq0\}=\{(\frac{x_0}{x_1},\ldots,\frac{x_n}{x_1})\}$. $\tilde\gamma$ is the line $(t,0,\ldots,0)$ on this chart. Note that this is a different parameterization of the line $\tilde\gamma$ from example \ref{tildeline1}. 
\end{ex}

\begin{ex}\label{tildenode1}


If $X\in\Pp^2$ is the nodal cubic 
$$
\begin{array}{rcl}
X &=& \{x_0x_2^2=x_1^3+x_1^2x_0\} \\
  &=& \{[u^3:t^2u-u^3:t^3-tu^2]\}
\end{array}
$$
and choose 
$$
\gamma=(x_0=1,\,x_1=\eps^2-1,\,x_2=\eps^3-\eps),
$$
a 3-arc over the smooth point $[1:-1:0]$.  $\tilde\gamma$ is $X$.  

On the other hand, choosing arcs over the singular point gives a slightly different picture.  The two arcs 
$$
\zeta=(x_0=1,\,x_1=\eps^2-2\eps,\,x_2=\eps^3-3\eps^2+2\eps),
$$ 
and 
$$
\zeta'=(x_0=1,\,x_1=\eps^2+2\eps,\,x_2=\eps^3+3\eps^2+2\eps),
$$ 
over the point $[1:0:0]$ corresponding to the two branch directions are distinct arcs in $\widehat X_3$.  There is no local reparameterization taking one to the other.  But the curves $\tilde\gamma$ and $\tilde\eta$ are both $X$.  In general, if a curve has an ordinary $n$--fold point, there will be $n$ arcs giving that curve via the tilde construction.  Arcs are local objects, and retain some information about the parameterization.
\end{ex}

For tangents, this situation is generalized by the following theorem, proved in \cite{H-M} and also in \cite{K-K} by more algebraic means, but needing some extra hypotheses on the family of curves.

\begin{thm}{\rm \cite{H-M, K-K}}\label{fanorattan} Suppose that $X$ is a complex projective manifold that is covered by rational curves (i.e. a Fano variety), and $H$ a maximal dominating family of minimal degree rational curves.
Write $H_x$ for the sub-family of curves passing through a point $x\in X$ and $\widetilde {H_x}$ for the normalization of such.  Then for general $x$ the tangent map 
$$ 
\tau_x : \widetilde{H_x} \dasharrow \widehat {X_{1,x}}, \quad \ell \mapsto \widehat{\ell_{1,x}} 
$$
is birational.
\end{thm}

Using arcs to replace tangents, one wishes to follow this example and extract higher order information about rational curves in a space $X$. 
It is natural to ask in what generality a result such as \ref{fanorattan} is true.  What properties of a variety are reflected in the images of the higher order Taylor series maps?  If the birational map between $H_x$ and $\widehat X_{1,x}$ is replaced by rational map, finite map, or just morphism, for what class of varieties does the statement hold?  

The tilde construction for projective space should be viewed as a prototype for a similar construction in a more general context.  The curve $\tilde\gamma$ is the ``simplest'' curve that is approximated by $\gamma$.  For arcs that are rationally generated, one expects to find at most a finite set of rational curves to chose from when looking for a ``simplest'' curve approximated by a given arc.  The proper generalization of \ref{fanorattan} would say that there is a unique such curve.  

Theorem \ref{uni_rga_rc} should be viewed as a step in this direction. It follows almost immediately from the definitions of the relevant properties.

\begin{defn} A variety $X$ is called {\em unirational} if there exists a dominant rational map $\psi:\Pp^n \dasharrow X$.   
\end{defn}  

Unirational varieties share many properties with rational varieties, but this class is significantly larger.  For example, Clemens and Griffiths showed that a general cubic threefold is not rational \cite{C-G}.  On the other hand, Koll\'ar proved that any cubic hypersurface is unirational \cite{Koll}.  

\begin{defn} A variety is {\em rationally connected } if any two points $p,q\in X$ can be joined by a rational curve. 
\end{defn}  

This notion was introduced by Koll\'ar-Miyaoka-Mori \cite{KMM},  and by Campana \cite{Camp}.  Rational connectedness is a useful property in that there are a number of good local descriptions.  One very useful tool when working with rational connectedness is the maximally rationally connected fibration.

\begin{defn}\cite[IV.5.3]{Koll} The {\em maximally rationally connected (MRC) fibration} is a rational map $\fee:X\dasharrow Q$ such that all fibers of $\fee$ are rationally connected, and there are no rational curves through a general point of $Q$. 
\end{defn}

\begin{thm}{\rm \cite[IV.5.4]{Koll}} If $X$ is a smooth proper variety then the MRC fibration exists and is unique up to birational equivalence.
\end{thm}

The proof of Theorem \ref{uni_rga_rc} in fact provides a slightly more refined result than the one stated in the introduction.

\begin{thm}
If $X$ is any variety, then
$$
{\rm \begin{tabular}{c} $X$ is \\unirational\end{tabular}}\implies{\rm \begin{tabular}{c} There exists $e$ such that \\for all $m$, $X$ has generically\\ $em$-rationally generated $m$-arcs\end{tabular}}
$$
and
$${\rm \begin{tabular}{c} $X$ has generically\\ rationally generated arcs\end{tabular}}\implies{\rm \begin{tabular}{c} $X$ is rationally\\ connected\end{tabular}}
$$ 
\end{thm}

\begin{pf} For the first statement, the idea of the proof is to pull an arc back to $\Pp^n$, extend to a curve there, and then push forward that curve.  

Choose a finite dominant rational map $\psi:\Pp^n\to X$.  There is an open set $U\subset X$ where $\psi:\psi^{-1}U\to U$ is \'etale.  Hence by \cite[Proposition 2.2]{Blik}, $(\psi^{-1}U)_m=U_m\times_U\psi^{-1}U$.  Thus, given a generic arc $\gamma\in \widehat X_m$, there exists an arc $\eta \in  \widehat \Pp^n_m$ such that  $\psi_m(\eta)=\gamma$.  

There is a degree $m$ rational curve $\tilde\eta\subset\Pp^n$ generating $\eta$.  So $\psi(\tilde\eta)$ is a rational curve on $X$ generating $\gamma$.  Setting $e=\deg(\psi_*\cO(1))$, $\deg\psi(\tilde\eta)=em$.  This proves the first implication.

To show that generic rational generation of arcs implies rational connectedness, one only needs that a generic arc is rationally generated for some $d$, but there need not be any global bound on $d$, even for fixed $m$.

Because any variety is birational to an open subset of a smooth proper variety, we may assume that $X$ is smooth and proper.

Consider the MRC fibration $\fee$ and suppose that $X$ is not rationally connected, so that this is not the trivial fibration $\fee:X\to pt$.  Suppose $\gamma$ is an arc that is not contained in the fibers of $\fee$.  That is, $\gamma\not\in F$ for all subvarieties contracted by $\fee$.  If $\gamma$ is a rationally generated arc, there exists a preimage of $\gamma$ via $\tau$ which is a rational curve that is not contained in a fiber of $\fee$.  Of course, this is contrary to the definition of $\fee$.  So $\fee$ must be the trivial fibration. \finpf
\end{pf}

As mentioned in the introduction, it is conjectured that unirationality and rational connectedness are distinct, however, there are no known examples to distinguish the two.

\section{The line locus}\label{locusproof}\label{linelocus} \label{intersect}

In this section theorem \ref{locus} is proved.  Essentially, because the arc spaces can be described so explicitly, when an $m$-arc is generated by a rational curve of degree $d$ it can be seen in the intersection ring of $\widehat X_m$.  For $\widehat X_1$, this is easy to calculate, but the intersection rings for higher arc spaces are
not known. 

Because of the isomorphism $X_1\iso TX$, the intersection ring $A(\widehat\Pp^n_1)$ is described by elementary intersection theory.
$$\begin{array}{rcl}
A(\widehat\Pp^n_1)
&=& A(\Pp^n)[j_1]/
\left(\sum_{i=0}^nc_i((T\Pp^n))j_1^{n-i}\right)\\
&&\\
           &=&\Z[j_0,j_1]/(j_0^{n+1},\sum_{i=0}^n\binom{n+1}{i}j_0^ij_1^{n+1-i}).
\end{array}
$$

\begin{claim}\label{lineclass}
A line $\ell\subset\Pp^n$ induces the cycle class
$$
[\hat\ell_1]=j_0^{n-1}(j_0+j_1)^{n-1}\in\widehat\Pp^n_1.
$$
Conversely, the cycle $Nj_0^{n-1}(j_0+j_1)^{n-1}$ is represented by the arcs of $N$ lines in $\Pp^n$.
\end{claim}

\begin{pf}
By \ref{jetideal}, the ideal of $\hat\ell_1$ is obtained from the
ideal of $\ell$ by extending the base to $\C[\eps]/(\eps^2)$ from $\C$.  The
$n-1$ generators of $I(\ell)$ are linear expressions.  Under the
extension, these become $n-1$ linear expressions in the weight 0 coordinates
and $n-1$ linear expressions in the weight 1 coordinates.

Write $h_i=[\{x_i=0\}]$ for the class associated to the vanishing of a
weight $i$ coordinate.  Then $[\hat\ell_1]=h_0^{n-1}h_1^{n-1}$.

Computing on a chart where
$z_0\neq0$, $h_0=[\{x_0=0\}]$, so $j_0=h_0$.  A divisor $D$ with $[D]=j_1$ is
given by the equation
$$\begin{array}{rcl}
D&=&\left\{\mbox{$t$ coefficient of } \frac {x(t)}{z(t)} = 0\right\}\\
&&\\
&=&\left\{\frac{x_1z_0-z_1x_0}{z_0^2} = 0 \right\}.\end{array}
$$ Hence $h_1+h_0-2h_0=j_1$, or 
$$
h_1=j_0+j_1.
$$

Conversely, by intersecting the class $Nh_00^{n-1}h_1^{n-1}$ with $h_0$, one sees that this is supported over a curve in $\Pp^n$ of degree $N$, and intersecting with the class $h_1$, that the arcs of this curve do not intersect a generic hyperplane in a fiber, and thus form a finite set.  This is only possible for the arcs generating a union of $N$ lines. 
\finpf\end{pf}

Observe that if two lines $\ell,\ell'\subset\Pp^n$ are distinct, then 
$$
\hat\ell_1\cap\hat\ell'_1=\emptyset
$$ since at any point where the lines intersect, the arcs are different. 

This fact is generally true in the sense that if $C$ and $C'$ are
distinct rational curves of degree at most $m$, then
$$
\widehat C_m\cap\widehat C'_m=\emptyset\subset\Pp^n_m.
$$

In the notation introduced in \ref{lineclass}, a presentation of the intersection ring $A(\widehat\Pp^n_1)$ more suited to the present context is
$$ 
A(\widehat\Pp^n_1) =
\Z[h_0,h_1]/(h_0^{n+1},h_1^n+h_1^{n-1}h_0+\cdots+h_1h_0^{n-1}+h_0^n)
$$

The calculations using this context are best illustrated by some familiar examples.  The application of this machinery is essentially the same in all cases.

\begin{ex}\label{cubic}


The generic cubic surface $X\subset\Pp^3$ has 27 lines, first enumerated in the correspondence of Cayley and Salmon in the 1850s.  To see them, write $X=\{f=0\}$ for a generic surface, where $f$ is a homogeneous degree 3 form in variables $x_0,x_1,x_2,x_3$. 
$$
f=\sum_{ a+b+c+d=3 } k_{abcd} x_0^ax_1^bx_2^cx_3^d
$$
Consider an arc $\gamma\in\widehat\Pp^3_1$, and the associated line $\tilde\gamma$ parameterized by $\theta(\gamma(t))$. 

For the line $\tilde\gamma$ to lie on $X$, $\theta(\gamma(t))$ must satisfy the appropriate equation: 
$$
\tilde f = \sum_{ a+b+c+d=3 } k_{abcd} (x_{00}+x_{01}t)^a (x_{10}+x_{11}t)^b (x_{20}+x_{21}t)^c (x_{30}+x_{31}t)^d = 0.
$$
Of course, a similar equation is used to determine when $\widehat{\tilde\gamma}_m\in\widehat X_m$, namely
$$
\sum_{ a+b+c+d=3 } k_{abcd} (x_{00}+x_{01}\eps)^a (x_{10}+x_{11}\eps)^b (x_{20}+x_{21}\eps)^c (x_{30}+x_{31}\eps)^d = 0.
$$

The existence of such a $\gamma'\in X_m$ is equivalent to the curve $\tilde\gamma$ contacting $X$ to order $m$.

The next step is to define polynomials $f^{(r)}$ by rewriting $\tilde f = 0$ as a polynomial in $t$ and collecting coefficients:
$$
\tilde f = \sum_{r=0}^\infty f^{(r)}t^r=0
$$
For this equation to be true for general $t$ (and hence $\tilde\gamma\in X$), each coefficient $f^{(r)}$ must vanish identically.  Observe that as $m$ increases, the number of conditions on the $x_{pq}$ given as nontrivial coefficients $f^{(r)}$ of this equation increases only until $m=3$, beyond which point all the coefficients are automatically zero and the list does not change.  In other words, once a line contacts a cubic surface to order 3, it lies on the surface, and contacts to all orders.

The upshot is that if $\widehat{\tilde\gamma}_3\in X_3$, then $\tilde\gamma\subset X$.  The locus $L$ where this is true is the intersection of the four divisors in $\widehat\Pp^3_1$ described by $f^{(0)}=0$, $f^{(1)}=0$, $f^{(2)}=0$ and $f^{(3)}=0$.  The point now becomes to calculate this intersection.

Write $L\subset X$ for the locus of lines.
$$
\begin{array}{rcl}
[L]&=&3h_0(2h_0+h_1)(h_0+2h_1)3h_1\\
        &=&18h_0^3h_1+45h_0^2h_1^2+18h_0h_1^3,\\
        &=&18h_0^3h_1+45h_0^2h_1^2+18h_0(-h_0h_1^2-h_0^2h_1-h_0^3),\\
        &=&27h_0^2h_1^2,\quad\mbox{\rm the 27 lines on the cubic, by \ref{lineclass}}.
\end{array}$$
\end{ex}

\begin{ex}\label{quintic}
For the quintic threefold of string theory fame \cite{CDGP}, the procedure is the same, replacing the number 3 with the number 5.  The divisors in question are
$$
(5-i)h_0+ih_1 \:\:\: i=0\ldots5
$$ 
and the intersection is represented by the cycle
$$
[L]=2875h_0^3h_1^3
$$ 
indicating the lines.  I do not know the first appearance of this number.  Saunders Mac Lane attributes it without reference to Harris in 1975, and it was certainly known to Clemens and Katz in the early 1980s. 

In this example, it is especially important to maintain the genericity assumption about the threefold.  For example, the Fermat quintic,
$$
F=\{x_0^5+x_1^5+x_2^5+x_3^5+x_4^5=0\}
$$
contains an infinite number of lines.  See \cite{Katz1}.

In addition to determining the lines that lie on the hypersurface itself, this methodology can be used to detect lines in the ambient space that intersect the hypersurface to some fixed degree.  For example, in \cite{Paci}, it is determined that the surface of lines that intersect the quintic to degree five (or are tangent to degree four) form a surface of degree 650.  This calculation can be done using arcs.

As above, an arc $\gamma\in\widehat\Pp^4_1$ determines a line $\tilde\gamma\subset\Pp^4$, which then has arc spaces of its own, $\widehat{\tilde\gamma}_m$.  When $\widehat{\tilde\gamma}_m\subset\widehat X_m$, then $\tilde\gamma$ is tangent to $X$ to degree $m$.  So to recover Pacienza's result, it suffices to find the class $[M]$ of arcs $\gamma\in\widehat\Pp^4_1$ where this inclusion holds.  

This is found by computing the product
$$
\prod_{i=0}^4(5-i)h_0+ih_1=650h_1^3h_0^2+1225h_1^2h_0^3+650h_1h_0^4,
$$
and then intersecting with $h_0^2$ to get the degree of the surface.
\end{ex}

\begin{ex}\label{delpezzo}
The del Pezzo surface of degree four is a complete intersection of type (2,2) in $\Pp^4$.  It is also amenable to the same methods.  In this case there are two defining equations, say $f$ and $g$, both of order two.  The divisor classes of $f^{(i)}$ and $g^{(i)}$ are the same, so the line locus is the product of two copies of each of the three divisor classes: $2h_0$, $h_1+h_0$ and $2h_1$, corresponding to $f^{(0)}$ and $g^{(0)}$, $f^{(1)}$ and $g^{(1)}$, and  $f^{(2)}$ and $g^{(2)}$, respectively.  The product is computed in the intersection ring $A(\widehat\Pp^4_1)$ 
$$
[L]=\left(\prod_{i=0}^2((2-i)h_0+ih_1)\right)^2.
$$ 
This comes out to $16h_0^3h_1^3$ indicating the 16 lines.
\end{ex}

With these examples in mind, we continue with the proof of \ref{locus}.

\begin{pf}[Theorem \ref{locus}] Consider a generic complete intersection $X$ in $\Pp^n$ and an arc $\gamma\in\widehat X_d$.  There is a unique curve $\tilde\gamma\subset\Pp^n$ as per the construction in \ref{tilde}.  We know that 
$$
\tilde\gamma\subset X \iff f(\theta(\gamma(t)))=0
$$ 
for all $f\in I(X)$ and $t\in\C$.  As in example \ref{cubic}, if
$$
f=\sum_\alpha k_{s(\alpha)} x_\alpha^{s(\alpha)}
$$
define 
$$
\tilde f=\sum_\alpha k_{s(\alpha)} ( \sum_{i=0}^d x_{\alpha i} t^i)^{s(\alpha)}
$$
and define $f^{(i)}$ via rewriting $\tilde f$:
$$
\tilde f=\sum_{i=0}^\infty f^{(i)}t^i.
$$

With this setup, the following equivalences can be spelled out explicitly
$$\begin{array}{rcl}
\tilde\gamma\subset X & \iff & f(\theta(\gamma(t))=0 \\
                      & \iff & f^{(i)}=0 \quad \forall i \\
                      & \iff & f^{(i)}=0 \quad i=0,\ldots,md \\
                      & \iff & \widehat{\tilde\gamma}_i\subset\widehat X_i \quad  i=0,\ldots,md
\end{array}$$ 

At this point it is key that $[f^{(i)}=0]\in A(\widehat\Pp^n_m)$ can be described explicitly, so set $m=1$.  For $f$ homogeneous of degree $d$, the $f^{(i)}$ are composed of monomials all having a total of $d$ variables and having weight $i$ when the weight of $x_{ji}$ is $i$.
$$
\begin{array}{c|c|c}
\mbox{\rm equation} & \mbox{\rm typical term} & \mbox{image in $A(\widehat\Pp^n_1)$}\\
\hline
f^{(0)}=0 & x_{j0}^d & dh_0 \\
f^{(1)}=0 & x_{j0}^{d-1}x_{j1} & (d-1)h_0 + h_1 \\
f^{(2)}=0 & x_{j0}^{d-2}x_{j1}^2 & (d-2)h_0 + 2h_1\\
\vdots    & \vdots & \vdots \\
f^{(d-1)}=0 & x_{j0}x_{j1}^{d-1} & h_0 + (d-1)h_1\\
f^{(d)}=0 & x_{j1}^d & dh_1
\end{array}
$$

To find the locus $L$ where the arcs generate lines on $X$, compute the product of all the classes in the left hand column coming from a set of generators of $I(X)$.  Since $X$ is a complete intersection of type $(d_1,\ldots,d_r)$, 
$$
I(X)=\langle f_1,\ldots,f_r \rangle
$$
with $\deg f_i=d_i$.  Thus 
$$
L  =  \displaystyle\prod_{j=1}^r\prod_{i=0}^{d_j}(d_j-i)h_0+ih_1
$$
as claimed.

If $\sum d_i=2(n-1)-r$ then the degree of $L$ is $2(n-1)$, in other words, $L\in A^{2(n-1)}(\widehat\Pp^n_1)$, which is
generated by $h_0^{n-1}h_1^{n-1}$ and $h_0^nh_1^{n-2}$.  However, the second term does not appear in $L$.  Looking at the calculation, note that because the expression for $[L]$ is symmetric in $h_0$ and $h_1$, previous to the reduction by $\sum_{i=0}^n h_0^ih_1^{n-i}=0$, the coefficient on $h_0^nh_1^{n-2}$ is the same as that on $h_0^{n-2}h_1^n$.  Applying that reduction introduces one term that cancels the $h_0^nh_1^{n-2}$ term and one term that modifies the $h_0^{n-1}h_1^{n-1}$ term.  The rest of the terms introduced by this substitution are zero because $h_0^{n+1}=0$.
$$\begin{array}{rcl}
[L] &=& \alpha h_0^nh_1^{n-2}+\beta h_0^{n-1}h_1^{n-1}+\alpha h_0^{n-2}h_1^n\\ 
    &=& \alpha h_0^nh_1^{n-2}+\beta h_0^{n-1}h_1^{n-1}+\alpha h_0^{n-2}(-\sum_{i=1}^nh_0^ih_1^{n-i})\\
    &=& (\beta-\alpha) h_0^{n-1}h_1^{n-1}
\end{array}
$$
This is the cycle class of $\beta-\alpha$ lines in the base space, and over each point in the base lying on one of these lines, an arc generating the line.  There is unfortunately no closed form known for $\beta-\alpha$ in general, even for the hypersurface case.
\finpf \end{pf}

The assumption of generality is needed because over degenerate points of the moduli space of complete intersections the continuity of the cycle class $L$ is not assured.

\section{Lines through a generic point}\label{cigenpt}

This part is motivated by the work of Landsberg.  In \cite{Land}, he proves

\begin{thm}{\rm \cite{Land}}\label{landsberg}
The number of lines, if finite, through a general point of a variety
of dimension $m$ is at most $m!$.  Also, this bound is sharp in the
case of a generic hypersurface $X$ of degree equal to $\dim(X)$ in
projective space $\Pp^{\dim X+1}$.
\end{thm}

The method of this paper can be used to calculate this number for general complete intersections.  Analysis of how the main intersection theory calculation varies in families, as can be found in \cite{Fult}, recovers Landsberg's result, and provides an improved sequence of bounds based on possession of finer information about the varieties in question. 

Determination of the number of lines through a generic point on a
complete intersection in $\Pp^n$ is an easy step once the locus of
lines is known.  Suppose $X$ is a complete intersection of dimension $m=n-r$ and degree $d$ with lines generated by the cycle of arcs $L\in A(\widehat\Pp^n_1)$. To look at lines through a generic point
$x\in X$, intersect with $m$ hyperplanes to get arcs living over $d$ points of $X$ and divide the resulting class by $d$.  
$$
L_x=Lh_0^m/d
$$

\begin{ex} A quadric surface has lines generated by $L=2h_0(h_0+h_1)2h_1$.  It is a degree 2 surface in $\Pp^3$, so the lines through a generic point are given by 
$$
\begin{array}{rcl}
L_x &=& Lh_0^2/2 \\
    &=& 2h_0^3h_1^2,
\end{array}
$$
which is what is expected.
\end{ex}

\begin{ex}\label{34in8}
The intersection of a cubic and a quartic in $\Pp^8$ generally has a finite number of lines through a general point.  This number is given by the following cycle class in $A(\widehat\Pp^8_1)$.
$$
\frac{h_0^6}{12}\prod_{i=0}^3(ih_0+(3-i)h_1)\prod_{i=0}^4(ih_0+(4-i)h_1)
= 144h_0^8h_1^7
$$
\end{ex}

A degree count in the general expression of $L$ shows that for a general complete intersection of type $(d_1,\ldots,d_r)$, the number of lines through a general point is finite only if $\sum_{i=1}^rd_i=n-1$ (such varieties are maximal degree Fano varieties), and in this case it is equal to
$\prod_{i=1}^rd_i!$.  This is interesting.  It gives a list of the
values that the number of lines through a general point of a general complete
intersection can take, and will lead to the improvements alluded to above.  

This count is not valid in all cases.  It can happen, such as in the case of a cubic threefold (which has 6 lines through a general point) degenerating to a union of three hyperplanes, that the dimension of the locus of lines through a general point may jump.  For the cubic threefold, this is not a problem, since we are interested at the outset only in varieties where the dimension of the line locus is the same as the dimension of the variety, and hence the number of lines through a generic point is finite.  However, as the next example shows, the dimension of this locus can jump from a number less than the dimension of the variety to equal to the dimension of the variety, giving varieties that are naively expected to have no lines through a general point, but in fact have a finite number.

\begin{ex}\label{quartic_lines}Consider $S\subset Q$, the surface of lines in a quartic three-fold.  One calculates that $\deg S=320$, and by Lefschetz's theorem (see \cite[Exercise III.11.6(c)]{Hart}) $S$ is a complete intersection.  However, there is one line through a generic point of $S$, namely the line in $Q$ through that point of $S$.
\end{ex}

However, such cases are the exception, and it is probable that with some careful analysis of the moduli space of complete intersections, that similar counts or bounds could be extended to these spaces.  For most complete intersections, the lines through a general point are counted by arcs, and fall into this surprisingly simple pattern.

\begin{ex}\label{sixfold}
If $X$ is a generic complete intersection sixfold with a finite number of lines through a generic point, then either

\begin{itemize}
\item $X$ is a sextic in $\Pp^7$ and has $6!=720$ lines through a
general point,
\item $X$ is the intersection of a quadric and a quintic in $\Pp^8$ and
has $2!5!=240$ lines through a general point,
\item $X$ is as in example \ref{34in8} and has $3!4!=144$ lines
through a general point,
\item $X$ is of type $(2,2,4)$ in $\Pp^9$ and has $2!2!4!=96$ lines
through a general point,
\item $X$ is of type $(2,3,3)$ in $\Pp^9$ and has $2!3!3!=72$ lines
through a general point,
\item $X$ is of type $(2,2,2,3)$ in $\Pp^{10}$ and has $2!2!2!3!=48$
lines through a general point, or
\item $X$ is of type $(2,2,2,2,2)$ in $\Pp^{11}$ and has
$2!2!2!2!2!=32$ lines through a general point.
\end{itemize}
\end{ex}

This calculation is also relevant to the study of varieties that are
not complete intersections.  Because any variety can be found as a
component of some complete intersection, this sequence of calculations
leads to the sequence of bounds in Theorem \ref{genpt}.

\begin{pf}[Theorem \ref{genpt}]
The argument proceeds in the following manner:  first the case where $X$ is a complete intersection is worked out.  Using theorem \ref{locus} the line locus is calculated, and as in the examples above, a further calculation gives the lines through a generic point.  When $X$ is not a complete intersection, at a smooth point (and a generic point is of course smooth), it is still a local complete intersection.  So $X$ can be included into a reducible complete intersection.  By doing the intersections in families, the explicit calculation that was done for generic complete intersections provides the upper bound as described.

When $X$ is a complete intersection, the result is simply a rephrasing of theorem \ref{locus}.  The line locus of $X$ is given by the cycle
$$ 
L=\prod_{j=0}^r\prod_{i=0}^{d_j}\left( ih_0+(d_j-i)h_1\right),
$$ 
and the degree is $d=\deg X=\prod_{i=1}^rd_i$.  A generic point on a degree $d$ dimension $m$ variety $X\subset\Pp^n$ is found cycle--theoretically by 
$$
[x]=\frac{[X]h_0^m}{d}.
$$
The lines through a generic point are $L_x=L\cdot[x]$: 
$$
\begin{array}{rcl}
L_x &=& \dfrac{\left({\displaystyle\prod_{j=1}^r\prod_{i=0}^{d_j}}\left( ih_0+(d_j-i)h_1\right)\right)h_0^{n-r}}{{\displaystyle\prod_{j=1}^r}d_j} \\
    &=& \dfrac{\left({\displaystyle\prod_{j=1}^r\sum_{k=0}^{d_j}} a_{jk}h_0^kh_1^{d_j-k}\right)h_0^{n-r}}{{\displaystyle\prod_{j=1}^r}d_j} \\
    &=& \dfrac{\left(\displaystyle\sum_{p=r}^{\sum d_j}b_ph_0^ph_1^{\sum(d_j+1)-p}\right)h_0^{n-r}}{\displaystyle\prod_{j=1}^rd_j}, \\
&=& \dfrac{b_r}{\displaystyle\prod_{j=1}^rd_j}h_0^nh_1^{n-1}.

\end{array}
$$

To calculate $b_r=\prod d_j!d_j$, observe that to get a monomial with only $h_0^r$, one must choose the $h_1$ term from every binomial $ih_0+(d_j-i)h_1$ in the product.

This takes care of the case of when $X$ is a complete intersection. 

An arbitrary variety $X$ is generically smooth, and at a smooth point, a variety is a local complete intersection \cite[II.8.17]{Hart}.  That is, if $x\in X$ is a general point, there exists an inclusion
$X\subset\widetilde X$ with $X$ a component of $\widetilde X$ and
$\widetilde X$ smooth at $x$.  

$\widetilde X$ is contained in the same cycle class as a smooth complete intersection $Y$ of the same type.  To fix notation, we consider a family
$$
\xymatrix{
\widetilde X  \ar[r] \ar[d] & \sY \ar[d] \\
\{0\} \ar[r] & T}
$$
and write $Y$ for a generic member of the family other than $X$.
The analysis now divides into three cases, depending on the comparative dimensions of the line locus $L_Y$ of $Y$ and $Y$ itself.

{\sl Case 1 ($\dim L_Y>\dim Y$):}  In this situation, there are an infinite number of lines through a general point of $Y$.  The deformation over $T$ also gives a deformation of the morphism $L_Y\to Y$.  The dimension of the fibers of this morphism are $\dim L_Y-\dim Y$, and cannot deform to zero dimensional fibers by semi-continuity \cite[III.12.8]{Hart}.  But if $\widetilde X$ has a finite number of lines through a general point, $\dim L_{\widetilde X}=\dim \widetilde X$.

{\sl Case 2 ($\dim L_Y=\dim Y$):}  This is the situation of principal interest.  It will be shown that the number of lines through a general point of $\widetilde X$ is either infinite or bounded above by the number of lines through a general point of $Y$. 

Fix $n,d_1,\ldots,d_r$ with $\sum_{i=1}^rd_i=n-1,\,d_i>1$.  The space
of all complete intersections of type $(d_1,\ldots,d_r)$ is a quotient $\sY=\Pp^{N_1}\times\cdots\times\Pp^{N_r}/\sim$; the equivalence relation contains the information whereby multiple sets of equations can describe the same variety (or algebraically, an ideal can have many possible generating sets).  When $X$ is a hypersurface, it is just a quotient by $\C^*$, but it increases in complexity as $r$ increases.

Write $\sX\subset\sY\times\Pp^n$ for the universal complete
intersection of type $(d_1,\ldots,d_r)$.  That is, for $p\in\sY$, the
fiber $X_p\subset\Pp^n$ is the complete intersection parameterized by
$p$.  Consider the incidence variety $\sV\subset\sX\times\widehat\Pp^n_1$ where $(p,\gamma)\in\sV$ when $\gamma(0)=x\in X_p$ and $\tilde\gamma$ is a line in $X_p$
through $x$. The key is that the projection $\sV\to\sX$ is generically
finite (So conservation of number applies), and the fiber over a point $(p,x)$ of this projection
describes the arcs extending to lines in $X_p$ through $x$.

For generic $p\in\sY$, $X_p$ is smooth, and the line locus computation gives the cycle class of the fiber $V_p$ of $\sV\to\sX$.  A
generic member of this cycle class is a finite set of points in $\sV$,
so the map is generically finite.  That is, there are a finite number of 1--arcs rooted at a particular point that extend to lines contained in a particular complete intersection.  

The number of points in a fiber of
a generically finite map can decrease (for example, two lines could
come together, and the map $\sV\to\sX$ become ramified), but the
number cannot increase without becoming infinite.  What this means is that as $p$ varies to give a deformation to some reducible $X_p$ where the variety under consideration is found as a component, the number of arcs extending to lines does not increase. 

{\sl Case 3 ($\dim L_Y<\dim Y$):}  This case requires more careful analysis, and is not included in the present result.  However, the generic behavior in this case is that there are no lines through a general point.  It is only in special cases such as the Fano surface of lines in a quartic threefold, \ref{quartic_lines}, where the intersections determining the line locus is not transverse, that there are any such lines.  
\finpf
\end{pf}

Applying this result to the example \ref{sixfold} above, if a generic sixfold in $\Pp^9$ has finitely many lines through a general point, then there are at most 96, and if it is not contained in an intersection of two quadrics, then there are at most 72.

\nocite{*} 
\bibliographystyle{skalpha} 
\bibliography{Treisman-Rationally_generated_arcs}

\end{document}

%% file: macros.tex
\newcommand{\Pp}{{\bf P}}
\newcommand{\A}{{\bf A}}
\newcommand{\C}{{\mathbb C}}
\newcommand{\Z}{{\mathbb Z}}

\newcommand{\cO}{{\mathcal O}}

\DeclareFontFamily{OMS}{rsfs}{\skewchar\font'60}
\DeclareFontShape{OMS}{rsfs}{m}{n}{<-5>rsfs5 <5-7>rsfs7 <7->rsfs10 }{}
\DeclareSymbolFont{rsfs}{OMS}{rsfs}{m}{n}
\DeclareSymbolFontAlphabet{\scr}{rsfs}

\newcommand{\sL}{\scr{L}}

\newcommand{\sV}{\scr{V}}

\newcommand{\sX}{\scr{X}}
\newcommand{\sY}{\scr{Y}}

\newcommand{\eps}{\varepsilon}
\newcommand{\fee}{\varphi}

\DeclareMathOperator{\spec}{Spec}
\DeclareMathOperator{\proj}{Proj}

\DeclareMathOperator{\Hom}{Hom}

\newcommand{\iso}{\cong}

\newcommand{\finpf}{\hspace{\fill}$\square$ \\}

\theoremstyle{definition}
\newtheorem{defn}{Definition}[section]
\theoremstyle{plain}
\newtheorem{thm}[defn]{Theorem}

\newtheorem{claim}[defn]{Claim}

\theoremstyle{definition}
\newtheorem{ex}[defn]{Example}
\theoremstyle{definition}
\newtheorem{notn}[defn]{Notation}
\theoremstyle{remark}